# STATELESS QUANTUM STRUCTURES AND EXTREMAL GRAPH THEORY

Václav Voráček

Center for Machine Perception, Department of Cybernetics,
Faculty of Electrical Engineering, Czech Technical University,
Technická 2, 166 27 Praha, Czech Republic
(e-mail: voracva1@fel.cvut.cz)

We study hypergraphs which represent finite quantum event structures. We contribute to results of graph theory, regarding bounds on the number of edges, given the number of vertices. We develop a missing one for 3-graphs of girth 4. As an application of the graph-theoretical approach to quantum structures, we show that the smallest orthoalgebra with an empty state space has 10 atoms. Optimized constructions of an orthomodular poset and an orthomodular lattice with no group-valued measures are given. We present also a handcrafted construction of an orthoalgebra with no group-valued measure; it is larger, but its properties can be verified without a computer.

## 1. Introduction

Shortly after the discovery of quantum mechanics, it was shown in [1, 2] that, unlike the classical case, the Boolean lattice is not sufficient for describing the logic of quantum mechanics. The suitable algebraic structures for the description of quantum logic of our interest (called *quantum structures*) are the following:

- orthomodular lattice (OML) [3, 4],
- orthomodular poset (OMP) [5],
- orthoalgebra (OA) [6, 7].

We continue the study of quantum structures and focus on measures on them. Quantum structures can be represented by hypergraphs. Measures (also called *states*) on quantum structures are solutions to systems of linear equations, with variables corresponding to vertices and equations corresponding to edges of its hypergraph representation. For the constructions of quantum structures with specific space of measures (e.g. no measures at all), it is appropriate to find quantum structures represented by hypergraphs with a high number of edges, given the number of vertices. We adopt this idea and answer the questions about how many edges such hypergraphs may have, given their number of vertices.

The paper is organized as follows. In Section 2, we introduce hypergraphs representing quantum structures. Then we show the results of extremal graph theory related to quantum structures. In Section 3, we show that the smallest known stateless OA is the smallest possible among those represented by Greechie diagrams. Then we show a handcrafted example of an OA with no group-valued measure. Section 4 shows our constructions of an OML and OMP with no group-valued measure. Finally, in Section 5, we summarize our results.

## 2. Hypergraphs

In this section, we define hypergraphs and some of their basic terms. We refer the reader for more details to [8]. Then we show graph-theoretical bounds on hypergraphs related to quantum structures.

This work summarizes and extends the results of the author's thesis [9].

### 2.1. Definitions

DEFINITION 1. A *hypergraph* is a pair $(V, E)$, where $V$ is a finite set of *vertices* and $E$ is a set of *edges* (nonempty subsets of $V$) such that $\bigcup_{e \in E} e = V$.

- A *cycle* of *order* $n \geq 2$ is a sequence of distinct edges and vertices $v_1, e_1, v_2, e_2, v_3, \ldots, v_n, e_n$ such that the following holds:
  * $\{v_1, v_2, \ldots, v_n\} \subseteq V$,
  * $\{e_1, e_2, \ldots, e_n\} \subseteq E$,
  * $v_i, v_{i+1} \in e_i$ for each $i \in \{1, 2, \cdots, n-1\}$ and $v_1, v_n \in e_n$.
- The *degree* of a vertex $v \in V$, written $d(v)$, is the number of edges $e \in E$ containing $v$.
- For a vertex $v \in V$, we denote by $\mathrm{adj}(v)$ the set of all vertices adjacent to $v$,
$$\mathrm{adj}(v) = \{u \in V \setminus \{v\} \mid \exists e \in E : \{u, v\} \subseteq e\}.$$

Following [10], finite quantum structures can be represented by hypergraphs. Their vertices are *atoms* of the quantum structure, i.e. nonzero elements $a$, which are minimal in the sense that there is no $b$ such that $0 < b < a$. The edges are sets of vertices corresponding to maximal sets of mutually orthogonal atoms. (Elements $a, b$ of a quantum structure are *orthogonal* iff $a \leq b'$; this relation is symmetric.) Edges correspond to maximal Boolean subalgebras of the quantum structure. Each finite quantum structure is uniquely determined by its hypergraph.

Not all hypergraphs represent quantum structures in the above manner. Then necessary and sufficient conditions for a hypergraph to represent a particular type of a quantum structure are rather complicated, see [4]. Here we reformulate simplified known conditions.

PROPOSITION 1. *The following condition* (1) *is necessary and the conjunction of conditions* (1) *and* (2) *is sufficient for a finite hypergraph* $H = (V, E)$ *to represent an orthoalgebra,*

$$\forall e, f \in E, \ e \neq f : |e \setminus f| \geq 2, \quad (1)$$

$$\forall e, f \in E, \ e \neq f : |e \cap f| \leq 1. \quad (2)$$

*When the condition* (1) *does not hold, some vertices cannot represent atoms. Condition* (2) *excludes cycles of order* 2; *in other words, it requires a hypergraph of girth* 3.

DEFINITION 2. Definition A finite hypergraph satisfying conditions (1) and (2) of Proposition 1 is called a *Greechie diagram*.

A Greechie diagram can contain two-element edges, but these must be disjoint to all other edges. We shall concentrate on Greechie diagrams such that each edge has at least three vertices. Then condition (1) follows from (2). Every Greechie diagram represents an OA.

PROPOSITION 2 ([10]). *A Greechie diagram represents*

- *an OMP iff it has no cycle of order less than* 4,
- *an OML iff it has no cycle of order less than* 5.

*In terms of graph theory, Greechie diagrams of girth* 4 (*resp.* 5) *represent OMPs* (*resp. OMLs*).

### 2.2. Measures

Quantum *states* on (finite) quantum structures are *probability measures*, i.e. additive functions that send the maximal element to 1. They are in a one-to-one correspondence with the following functions defined on the respective hypergraphs (here extended to general commutative groups).

DEFINITION 3. Let $H = (V, E)$ be a hypergraph.

- A mapping $m : V \to [0, 1]$ is a *probability measure* (*state*) on $H$, if for every edge $e \in E$,
$$\sum_{v \in e} m(v) = 1.$$

- Let $G$ be a commutative group. A mapping $m : V \to G$ is a *group-valued measure* on $H$, if there is an element $g \in G$ such that for each edge $e \in E$,
$$\sum_{v \in e} m(v) = g.$$

Both introduced measures are solutions to systems of linear equations. A group-valued measure can be seen as a solution to a *homogenous* system of linear equations with $|V|+1$ variables, $g$ and $m(v)$, $v \in V$, and $|E|$ equations (one for every edge). Clearly, every quantum structure admits the trivial group-valued measure, assigning 0 to every vertex; we are interested in other solutions. For the rest of the paper we will for brevity write "group-valued" instead of "nontrivial group-valued".

LEMMA 1 ([11]). *Let H be a hypergraph. If there is a group-valued measure on H, then there is a group-valued measure on H, whose range is a subset of the cyclic group $\mathbb{Z}_p$ for some prime p.*

This allows us to restrict our attention to cyclic groups $\mathbb{Z}_p$ for some prime $p$. Given a hypergraph, it is easy to determine whether there is a probability measure on it. The problem is more complicated for group-valued measures. If the number of edges equals the number of vertices plus one, then the matrix of the system of linear equations is a square matrix, and it is sufficient to check whether its determinant is $\pm 1$. We provide a slightly more general lemma.

LEMMA 2. *A system of linear equations $\mathbf{A}x = \mathbf{0}$, where $\mathbf{A}$ is an $m \times m$ square matrix with integer entries, has a nontrivial solution from some $\mathbb{Z}_p^m$ iff $|\det(\mathbf{A})| \neq 1$.*

*Proof:* Let us inspect several cases.

- $|\det(\mathbf{A})| = 1$.
  Then there exists $\mathbf{A}^{-1}$, its entries are quotients of minors of $\mathbf{A}$ and $\det(\mathbf{A})$, so they are integers. Multiplying a matrix by its inverse can be seen as applying the needed steps of Gauss–Jordan elimination to obtain the identity matrix. Since $\mathbf{A}^{-1}$ has only integer entries, the elimination was divisionless, and there is no nontrivial solution to the system of equations.
- $\det(\mathbf{A}) = 0$.
  Then the null-space of $\mathbf{A}$ is nontrivial, and there exists a nontrivial solution to the system of equations.
- $\det(\mathbf{A}) \notin \{-1, 0, 1\}$.
  Then $\det(\mathbf{A})$ is divisible by a prime $p$, and the matrix of the system has not the full rank when considered as a mapping of a vector space over the group $\mathbb{Z}_p$. Hence there is a nontrivial solution to the system of equations. □

If there are at most as many edges as vertices, then the second case of the previous proof still applies.

COROLLARY. *A hypergraph with no group-valued measures has more edges than vertices.*

### 2.3. Bounds

We show upper bounds on the number of edges, given the number of vertices for Greechie diagrams. If we add a two-element edge, it must be disjoint to all other edges, thus we added two vertices and one edge. This would not improve the bound, thus we may exclude two-element edges. On the other hand, edges with more than three vertices would unnecessarily increase the number of vertices. Looking for extremal hypergraphs, we may, without loss of generality, restrict attention to 3-graphs, i.e. assume that every edge contains exactly 3 vertices. From now on, when no ambiguity may occur, we will study a hypergraph $H = (V, E)$ and denote $m = |E|$ and $n = |V|$.

THEOREM 1. *Every Greechie diagram of an OA satisfies*

$$m \leq \frac{n(n-1)}{6}.$$

The proof can be found in [12]. The bound is tight and can be attained whenever $n(n-1)/6$ is an integer.

THEOREM 2. *Every Greechie diagram of an OMP satisfies*

$$m \leq \frac{n(n+3)}{18}.$$

*Proof:* Since $H$ is a 3-graph, it clearly holds that $|\text{adj}(v)| = 2d(v)$ for any vertex $v \in V$.

Now consider an edge $\{v_1, v_2, v_3\} \in E$. Then

$$|\text{adj}(v_1)| + |\text{adj}(v_2)| + |\text{adj}(v_3)| \leq n + 3.$$

This follows from the fact that there is no cycle of order 3 or less, therefore $\text{adj}(v_1) \cap \text{adj}(v_2) = \{v_3\}$, etc.

Hence

$$2d(v_1) + 2d(v_2) + 2d(v_3) \leq n + 3,$$

$$d(v_1) + d(v_2) + d(v_3) \leq \frac{n+3}{2}.$$

Summing up these inequalities for every edge, we obtain the following

$$\sum_{v \in V} d(v)^2 \leq m \frac{n+3}{2}.$$

Now apply the root mean square-arithmetic mean (RMS-AM) inequality

$$\frac{(\sum_{v \in V} d(v))^2}{n} \leq \sum_{v \in V} d(v)^2 \leq m \frac{n+3}{2}.$$

Note that $\sum_{v \in V} d(v) = 3m$, then we can write

$$\frac{9m^2}{n} \leq m \frac{n+3}{2},$$

$$m \leq \frac{n(n+3)}{18}. \qquad \square$$

The bound is not asymptotically optimal. It was shown in [13] that for any positive constant $c$, there exists a number of vertices $n$ such that $m \leq cn^2$. On the other hand, the best polynomial upper bound is quadratic, so the proposed bound is reasonable even for high values of $n$. We are interested in those values of $n$, for which $n$ is close to $m$; then the bound is exact, if the RMS-AM inequality may, in fact, be an equality. This may happen for $n \in \{3, 9, 15\}$.

THEOREM 3. *Every Greechie diagram of an OML satisfies*

$$m \leq \frac{2n\sqrt{n - \frac{3}{4}} + n}{12}.  \qquad (3)$$

The proof can be found in [14]. This bound is asymptotically correct, and is almost tight even for small values of $n$. According to this formula, the lowest value of $n$, for which the bound may be greater than or equal to $m$, is 31. In [15], there were generated all small Greechie diagrams for OMLs and it turned out that the number of vertices needed is, in fact, 35. For $n = 35$, inequality (3) tells us that $m \leq 38$.

## 3. Orthoalgebras

In [6, 16], there are examples of OAs admitting no probability measure, having 10 atoms. We show that these examples are minimal.

THEOREM 4. *A Greechie diagram $H = (V, E)$ of an OA with $|V| \leq 9$ always admits a probability measure.*

*Proof:* We will enumerate all possible OAs with 9 vertices and construct probability measures on them. The proof of the case when $|V| \leq 8$ clearly follows from the proof for $|V| = 9$. Therefore we assume that $|V| = 9$.

The upcoming diagrams will capture a probability measure. Every symbol for a vertex has assigned a value of a probability measure. Every edge is represented by a maximal line segment.

First, consider the case, where every edge has precisely 3 vertices. Then there is a probability measure assigning $1/3$ to every vertex.

There clearly cannot be 4 edges containing 4 or more vertices. Now inspect the remaining cases.

- There are three edges $e, f, g \in E$ with 4 vertices. Then there is the following probability measure,

  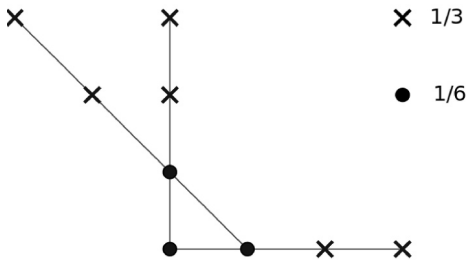

  If there are other edges in the hypergraph, then they contain three crosses, $1/3 + 1/3 + 1/3 = 1$.
- There are precisely two edges $e, f \in E$ with 4 or more vertices.

- If $e \cap f = \emptyset$, then there is at most one vertex $v$ not contained in these edges. Hence if there is any edge in $E$ distinct from $e$ and $f$, it contains precisely one vertex from $e$, one vertex from $f$, and $v$. There is a probability measure which assigns $1/|e|$ and $1/|f|$ to all vertices in $e$ and $f$, respectively. If $v$ exists, its measure is $1 - 1/|e| - 1/|f|$.
- Otherwise, $e \cap f = \{v\}, v \in V$.

    * If $|e|+|f| \geq 9$, then there is at most one vertex $u \in V$ not contained in $e$ and $f$. Hence there is no edge beside $e, f$ containing $v$. Therefore if there is another edge in $E$, then it consists of $u$, an element of $e \setminus \{v\}$, and an element of $f \setminus \{v\}$. Then there is an analogous probability measure to the previous case, which assigns 0 to $v$.
    * Otherwise, $|e| = |f| = 4$. There are two vertices $u, w \in V$ not contained in $e$ and $f$. If there is no other edge containing $v$, we can say the measure of $v$ is 0, and the measure of all other elements is $1/3$. Otherwise, there is an edge $\{u, v, w\} \in E$. Then there is the following probability measure,

    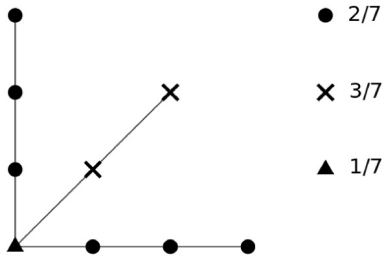

    If there are other edges in the hypergraph, they contain two circles and a cross, $2/7 + 2/7 + 3/7 = 1$.

- Let us assume that there is a unique edge $e \in E$ with 4 or more vertices.

    - If there is no edge $f \in E$ such that $e \cap f = \emptyset$, then there is a probability measure which assigns $1/|e|$ to all vertices contained in $e$ and $(1-1/|e|)/2$ to the others.
    - Otherwise, there is an edge $f \in E$ such that $e \cap f = \emptyset$. If there is precisely one such edge, then there are two more cases. The first case is

    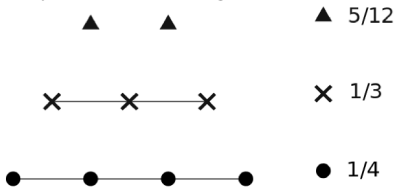

    If there are other edges in the hypergraph, then they contain a circle, a triangle (the case when an edge contains two triangles is treated in the

sequel) and a cross vertex, $1/4 + 1/3 + 5/12 = 1$. The second case is

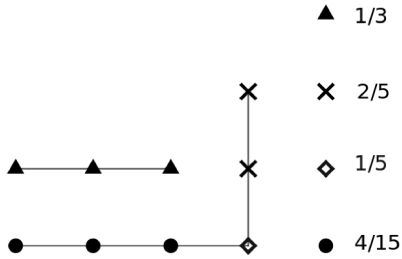

If there are other edges in the hypergraph, then they contain a triangle, a cross and a circle, $1/3 + 2/5 + 4/15 = 1$.

– If there are two edges with an empty intersection with $e$, then there exists the following probability measure,

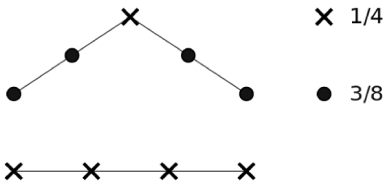

If there are other edges in the hypergraph, then they contain two circles and a cross, $3/8 + 3/8 + 1/4 = 1$. □

Every known example of a quantum structure with no group-valued measure (those in [11, 16, 17]) used a computer for generating the quantum structure, and/or for the verification that there is no group-valued measure, so the actual idea behind the construction remains hidden. We present a handcrafted example of an OA with no group-valued measure.

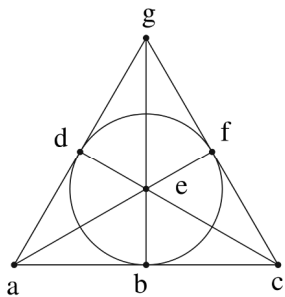

**Fig. 1.** An OA with 7 edges and 7 vertices.

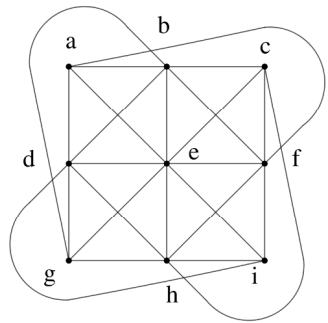

**Fig. 2.** An OA with 12 edges and 9 vertices.

LEMMA 3. *If there is a nonconstant group-valued measure on the OA captured in Fig. 1, then its range is isomorphic to* $\mathbb{Z}_2$.

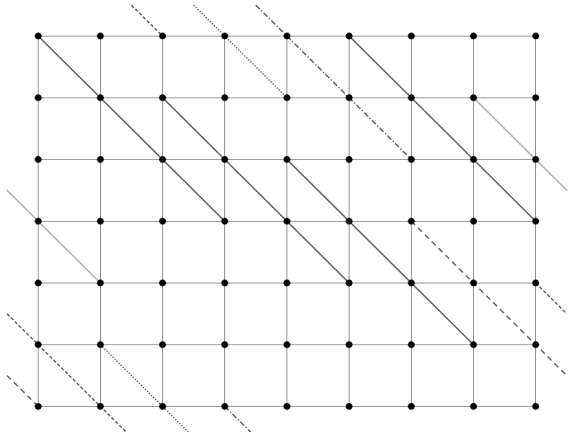

**Fig. 3.** An OA with no group-valued measures.

*Proof:* Take arbitrary two vertices, e.g. $g$ and $b$, and sum up measures over all edges containing them.

Let $S = \sum_{v \in V} m(v)$.

- $[m(g)+m(d)+m(a)]+[m(g)+m(e)+m(b)]+[m(g)+m(f)+m(c)] = S+2m(g)$,
- $[m(b)+m(a)+m(c)]+[m(b)+m(e)+m(g)]+[m(b)+m(d)+m(f)] = S+2m(b)$.

Hence $S + 2m(g) = S + 2m(b)$ and $2(m(g) - m(b)) = 0$, so either $2 = 0$ or $m(g) = m(b)$. The other cases are solved similarly. □

LEMMA 4. *If there is a nonconstant group-valued measure on the OA captured in Fig.* 2*, then its range is isomorphic to* $\mathbb{Z}_3$.

The proof is analogical to the previous case and therefore, is omitted.

Let us describe the diagram of an OA in Fig. 3. It is a simplified Greechie diagram. The vertices of a column form an OA from Fig. 1, and the vertices of a row form an OA from Fig. 2. Edges of these hypergraphs are truncated to vertical and horizontal lines, respectively. The diagonal edges are edges in the normal sense. For clarity we draw all the edges as straight lines. If an edge could not be drawn as a straight line, we drew it as multiple disconnected straight lines using the same and unique style.

THEOREM 5. *The OA from Fig.* 3 *does not admit any group-valued measure.*

*Proof:* Let there be a group-valued measure $m$ on it. If $m$ is constant, then it is trivial, as there are edges with both, 3 and 4 vertices. Hence $3m(v) = 4m(v)$ and $m(v) = 0$ for any vertex $v \in V$.

Assume $m$ is not constant, then $m$ is either $\mathbb{Z}_2$- or $\mathbb{Z}_3$-valued.

Let $m$ be $\mathbb{Z}_2$-valued, the other case is analogical and therefore omitted. Then the measures of vertices are constant in every row. Let $r_1$ be the measure of vertices in the first row, $r_2$ in the second row, etc.

There is a diagonal edge with vertices from 1st, 2nd, 3rd, and 4th row and another edge with vertices from 2nd, 3rd, and 4th and 5th, row. Hence

$$r_1 + r_2 + r_3 + r_4 = r_2 + r_3 + r_4 + r_5,$$
$$r_1 = r_5.$$

Similarly, we get $r_1 = r_5 = r_2 = r_6 = r_3 = r_7 = r_4$. Hence $m$ is a constant measure and the OA does not admit any group-valued measure. □

The technique of combining OAs as shown in Fig. 3 without the diagonal edges also work for OMPs.

REMARK 1. The construction which merges two OAs, as used in Fig. 3 can be also used to merge two OMPs.

## 4. Orthomodular posets and lattices

We optimized the constructions of OML and OMP with no group-valued measures. The former smallest known such OML had 74 vertices, here is a construction with just 67 vertices. The construction is still far from the graph-theoretical lower bound, so there is hope for improvement. For the OMP, there is still space for improvement of either construction, or lower bound. We conjecture that the smallest possible OMP with no group-valued measures has 20 or 21 vertices.

THEOREM 6. *The following quantum structures do not admit any group-valued measure.*

- *The OML having the Greechie diagram with vertices* $1, 2, \ldots, 67$ *and edges:*
  $\{1, 5, 65\}$, $\{1, 9, 21\}$, $\{1, 45, 56\}$, $\{2, 3, 18\}$, $\{2, 23, 64\}$, $\{2, 28, 40\}$, $\{3, 7, 17\}$, $\{3, 38, 63\}$, $\{4, 5, 6\}$, $\{4, 8, 24\}$, $\{4, 46, 62\}$, $\{5, 13, 37\}$, $\{5, 51, 52\}$, $\{6, 7, 60\}$, $\{6, 15, 31\}$, $\{6, 20, 35\}$, $\{7, 11, 32\}$, $\{7, 44, 67\}$, $\{8, 11, 12\}$, $\{8, 27, 29\}$, $\{9, 12, 18\}$, $\{9, 14, 20\}$, $\{10, 20, 28\}$, $\{10, 26, 32\}$, $\{10, 46, 63\}$, $\{11, 19, 25\}$, $\{12, 42, 49\}$, $\{13, 17, 22\}$, $\{13, 42, 58\}$, $\{14, 19, 58\}$, $\{14, 23, 36\}$, $\{15, 26, 29\}$, $\{15, 39, 43\}$, $\{16, 19, 65\}$, $\{16, 27, 35\}$, $\{16, 38, 39\}$, $\{17, 21, 29\}$, $\{19, 31, 59\}$, $\{21, 46, 53\}$, $\{22, 24, 28\}$, $\{22, 45, 48\}$, $\{23, 49, 52\}$, $\{24, 36, 38\}$, $\{25, 28, 50\}$, $\{25, 30, 47\}$, $\{27, 40, 57\}$, $\{30, 43, 61\}$, $\{31, 34, 48\}$, $\{32, 41, 52\}$, $\{33, 50, 54\}$, $\{33, 57, 63\}$, $\{33, 58, 60\}$, $\{34, 37, 44\}$, $\{34, 49, 50\}$, $\{35, 42, 47\}$, $\{39, 41, 54\}$, $\{40, 41, 59\}$, $\{40, 44, 56\}$, $\{42, 56, 62\}$, $\{43, 49, 55\}$, $\{45, 55, 60\}$, $\{46, 55, 59\}$, $\{47, 63, 66\}$, $\{48, 51, 53\}$, $\{51, 57, 61\}$, $\{53, 54, 67\}$, $\{61, 62, 64\}$, $\{64, 65, 66, 67\}$.
- *The OMP having the Greechie diagram with vertices* $1, 2, \ldots, 21$ *and edges:*
  $\{1, 4, 7\}$, $\{1, 5, 18\}$, $\{1, 10, 12\}$, $\{2, 3, 7\}$, $\{2, 6, 10\}$, $\{2, 17, 18\}$, $\{3, 5, 16\}$, $\{3, 8, 20\}$, $\{3, 14, 15\}$, $\{4, 6, 21\}$, $\{4, 8, 17\}$, $\{4, 11, 15\}$, $\{5, 6, 13\}$, $\{7, 13, 19\}$, $\{8, 9, 10\}$, $\{9, 11, 13\}$, $\{9, 16, 21\}$, $\{10, 15, 19\}$, $\{11, 12, 20\}$, $\{12, 16, 17\}$, $\{13, 14, 17\}$, $\{18, 19, 20, 21\}$.

*Proof:* The measure is a solution to a system of equations $\mathbf{Ax} = \mathbf{0}$. Both of the matrices $\mathbf{A}$ corresponding to these quantum structures have determinants $\pm 1$. According to Lemma 2, these quantum structures do not admit any group-valued measures. □

## 5. Summary

In the paper, we summarized results from extremal graph theory related to quantum structures and derived the missing one. We showed that the smallest known stateless OA is the smallest possible among a certain class of OAs. A handcrafted construction of an OA with no group-valued measures was given and also constructions of such an OMP and an OML.

We restricted attention to Greechie diagrams, thus excluding hypergraphs with cycles of order 2, or, equivalently, with edges intersecting in more than one vertex. Such hypergraphs also may represent quantum structures, even OMLs [4]. They are rather large and complicated, thus we conjecture that they do not exceed the bounds which we derived for Greechie diagrams. We proved this for OAs without states.

**Acknowledgement**

I thank prof. M. Navara for his help. This work was supported by the Czech Science Foundation grant 19-09967S.